\documentclass[reqno]{amsart}

\usepackage[all]{xy}

\numberwithin{equation}{section}

\theoremstyle{plain}
\newtheorem{lemma}{Lemma}[section]

\newtheorem{propo}[lemma]{Proposition}
\newtheorem{coro}[lemma]{Corollary}
\newtheorem{con}[lemma]{Conjecture}

\theoremstyle{definition}
\newtheorem{defi}[lemma]{Definition}

\newtheorem{example}[lemma]{Example}

\theoremstyle{remark}
\newtheorem{remark}[lemma]{Remark}
\newtheorem{notation}[lemma]{Notation}

\newcommand{\pic}{{\rm Pic\thinspace}}
\newcommand{\bs}{{\rm Bs\thinspace}}

\newcommand{\ac}{{\rm\bf A}}

\newcommand{\p}{\mathbb{P}}
\newcommand{\co}{\mathbb{C}}
\newcommand{\z}{\mathbb{Z}}

\newcommand{\oc}{{\mathcal O}}

\newcommand{\ls}{{\mathcal L}}

\newcommand{\ci}{{\mathcal I}}

\newcommand{\cre}{{\rm Cr\thinspace}}

\newcommand{\I}{{\mathcal I}}

\newcommand{\rt}{\longrightarrow}
\newcommand{\rmap}{\dashrightarrow}
\newcommand{\rr}{\rightsquigarrow}
\newcommand{\rri}{\stackrel{(i)}{\rr}}
\newcommand{\rrii}{\stackrel{(ii)}{\rr}}
\newcommand{\rriii}{\stackrel{(iii)}{\rr}}

\begin{document}

\title{On a class of special linear systems of $\p^3$}
\author{Antonio Laface}
\address{
Dipartimento di Matematica, Universit\`a degli Studi di Milano,
Via Saldini 50, \newline 20100 Milano, Italy }
\email{antonio.laface@unimi.it}

\author{Luca Ugaglia}
\address{
Dipartimento di Matematica, Universit\`a degli Studi di Milano,
Via Saldini 50, \newline 20100 Milano, Italy }
\email{luca.ugaglia@unimi.it}

\begin{abstract}
In this paper we deal with linear systems of $\p^3$ through fat
points. We consider the behavior of these systems under a
cubo-cubic Cremona transformation that allows us to produce a
class of special systems which we conjecture to be the only ones.
\end{abstract}
\maketitle

\section*{Introduction}

Let us take the projective space $\p^n$ and let us consider the
linear system of hypersurfaces of degree $d$ having some points of
fixed multiplicity. The virtual dimension of such systems is the
dimension of the system of degree $d$ polynomials minus the
conditions imposed by the multiple points and the expected
dimension is the maximum between the virtual one and $-1$. The
systems whose dimension is bigger than the expected one are called
{\em special systems}.\\
There exists a conjecture due to Hirschowitz (see \cite{hi}),
characterizing special linear systems on $\p^2$, which has been
proved in some special cases (see~\cite{cm,cm2,mi,ev}).
\\
In this paper we describe a class of special linear systems on
$\p^3$. The main tool will be the cubo-cubic Cremona
transformation~\ref{cubic} which allows us to transform a linear
system into another one. The dimension of the two systems is the
same, while the virtual one may be different. This is a new
phenomenon which does not occur in $\p^2$. In
Proposition~\ref{vir-change} we give a formula expressing the
difference between these virtual dimensions and in particular we
prove that if a transformation decreases the degree of a system
then it does not decrease its virtual dimension
(Corollary~\ref{decreasing-degree}). We will say that a system for
which is no longer possible to decrease the degree using one of
these transformations is in {\em standard form}. Starting from a
special system in standard form it is possible to construct
infinitely many special systems by applying a sequence of cubic
Cremona transformations. In this paper we describe two types of
special systems in standard form and we conjecture that they
produce all the possible special cases.
\\
The paper is organized as follows: in the first section we recall
some definitions and notations. In Section $2$ we give a
description of a cubic Cremona transformation of $\p^3$ and its
action on linear systems, while Section $3$ deals with the
resolution of the indeterminacy of this transformation. In Section
$4$ we state the conjecture and we give a procedure for evaluating
the dimension of a system and in the next one we give some
motivations for it. Section $6$ deals with special homogeneous
linear system according to the conjecture and finally Section $7$
provides some examples related to the procedure.

\section{Preliminaries}

We start by fixing some definitions and notations.

\begin{defi}
Given a sheaf $\mathcal{F}=\oc_{\p^3}(d)\otimes\I_V$, where $\I_V$
is the ideal sheaf of a subscheme $V\subset\p^3$, we denote by
$v(\mathcal{F})$ its {\em virtual dimension}, defined as
\[
v(\mathcal{F})=\chi(\mathcal{F})-1.
\]
\end{defi}

With $\ls=\ls_3(d,\allowbreak m_1,\allowbreak\ldots,\allowbreak
m_r)$ we will denote the linear system associated to the sheaf
$\oc_{\p^3}(d)\otimes\I_Z$, where $Z=\sum m_ip_i$ is a
zero-dimensional scheme of fat points. If no confusion arises, by
abuse of notation we will use the same letter $\ls$ to denote also
the sheaf. From the cohomology exact sequence associated to
\[
\xymatrix{ 0 \ar[r] & \ls\ar[r] & \oc_{\p^3}(d)\ar[r] & \oc_Z
\ar[r] & 0,}
\]
we obtain that $h^i(\ls) = 0$ for $i=2,3$. Therefore we deduce
that $v(\ls) = h^0(\oc_{\p^3}(d))-h^0(\oc_z)-1$, which may be
written also as
\begin{eqnarray*}\label{virt}
v(\ls) & = & \binom{d+3}{3}-\sum_{i=1}^r\binom{m_i+2}{3}-1.
\end{eqnarray*}

Let us denote by $e(\ls)=\max(v(\ls),-1)$ the {\em expected
dimension} of $\ls$.
\begin{defi}
A non-empty linear system $\ls$ is {\em special} if its expected
dimension is strictly smaller than the effective one or, which is
the same, if $h^1(\ls)\neq 0$.
\end{defi}
Let $(X,\pi)$ be the blow-up of $\p^3$ along $\{p_1,\ldots,p_r\}$;
by abuse of notation we will denote by $\ls$ the linear system
associated to $L=dH-\sum m_iE_i$, where $H$ is the pull-back of an
hyperplane of $\p^3$ and $E_i=\pi^{-1}(p_i)$. Let $\langle H, E_1,
\ldots, E_r\rangle$ and $\langle h, e_1,\ldots,e_r\rangle$ be two
bases for the Chow groups $A^1(X)$ and $A^2(X)$ respectively,
where $h$ is the pull-back of a line and $e_i$ is the class of a
line in $E_i$. The intersection matrix, with respect to these two
bases, is diagonal with the first element equal to $1$ and the
others equal
to $-1$ (since $E_i e_i=-1$).\\
Given a curve $C\subset\p^3$, by abuse of notation we will denote
by $\ls C$ the intersection product of their strict transforms in
$X$. We will write $C\in\ell_3(\delta,\mu_1,\ldots,\mu_r)$ to
denote a curve of degree $\delta$ with multiplicity $\mu_i$ at
$p_i$. In this way, the intersection is given by the formula
\begin{equation}\label{LdotC}
\ls C=d\delta-\sum_{i=1}^r \mu_i m_i.
\end{equation}
We recall the Riemann-Roch formula for a divisor $L$ on a smooth
threefold $X$:
\[
\chi(L) = \frac{L(L-K_X)(2L-K_X)+c_2(X) L}{12} + \chi(\oc_X).
\]

If the linear system associated to $L$ can be written as
$|L|=F+|M|$, where $F$ is a fixed divisor of $|L|$ and $|M|$ is
the residual system, then the above formula implies:
\begin{equation}\label{add}
v(L) = v(M) + v(F) + \frac{FM(L-K_X)}{2}.
\end{equation}

\section{Cubic Cremona transformations and linear systems}
In this section we focus our attention on a class of cubic Cremona
transformations of $\p^3$. Consider the system $\ls_3(3,2^4)$, by
putting the four double points in the fundamental ones, the
associated rational map is given by:
\begin{equation}\label{cubic}
\cre: (x_0:x_1:x_2:x_3) \rmap
(x_0^{-1}:x_1^{-1}:x_2^{-1}:x_3^{-1}).
\end{equation}

The birational map~\ref{cubic} induces an action on the pPicard
group of $X$ which can be described in the following way:
\begin{propo}\label{cre-surfaces}
The action of transformation~\ref{cubic} on
$\ls=\allowbreak\ls_3(d,\allowbreak m_1,\allowbreak \ldots,m_r)$
is given by:
\begin{eqnarray}\label{cre-a1}
\cre(\ls) & = & \ls_3(d+k,m_1+k,\ldots, m_4+k,m_5,\ldots,m_r),
\end{eqnarray}
where $k=2d-\sum_{i=1}^4 m_i$.
\end{propo}
\begin{proof}
Since the transformation~\ref{cubic} is an isomorphism out from the base locus,
we can reduce to the case $r=4$.
The set of monomials which generate $\ls$ (i.e. monomials of
degree $d$ in $x_1,\ldots,x_4$ with multiplicity $m_i$ at $p_i$
for $i=1\ldots,4$) can be represented by the set of their
exponents
\begin{equation}\label{T}
T:=\{(a_1,\ldots,a_4)\in\z^4\ |\ 0\leq a_i\leq d-m_i\
\textrm{and}\ \sum_{i=1}^4 a_i=d\}.
\end{equation}
In the same way $T'$ corresponds to $\ls_3(d+k,m_1+k,\ldots,
m_4+k)$. Observe that
\begin{eqnarray*}
\cre(\prod_{i=1}^4{x_i}^{a_i}) & = & \prod_{i=1}^4{x_i}^{d-a_i},
\end{eqnarray*}
where the right term may be written as
\[
(x_1^{m_1}\cdots x_4^{m_4})\prod_{i=1}^4{x_i}^{d-a_i-m_i}.
\]
This may be summarized by saying that the
transformation~\ref{cubic} induces the map $f(a_1,\ldots, a_4) =
(d-a_1-m_1,\ldots,d-a_4-m_4)$ from $\z^4$ to $\z^4$. We want to
prove that $f$ is a bijection between $T$ and $T'$. First of all
observe that if $(a_1,\ldots,a_4)\in T$ then $0\leq d-a_i-m_i\leq
(d+k) - (m_i+k)$ and $\sum_{i=1}^4 (d-a_i-m_i) = d + k$, which
implies that $f(T)\subseteq T'$. In the same way it is possible to
prove that $f(T')\subseteq T$ and, since $f^2=\textbf{1}_{\z^4}$,
we get the thesis.
\end{proof}

Observe that $\dim \cre(\ls) = \dim \ls$ but in general the
virtual dimensions of the two systems may be different.

\begin{lemma}\label{fixed-plane}
Suppose that $2d-m_1-m_2-m_3 < 0$ then $\bs\ls_3(d,m_1,m_2,m_3)$
contains the plane through the three points.
\end{lemma}
\begin{proof}
Consider $T$ as before (with $m_4=0$), then we have the inequality
$a_1+a_2+a_3 \leq 3d - m_1-m_2-m_3 < d$. This implies that each
monomial contains $x_4$ which gives the fixed plane.
\end{proof}

\begin{propo}\label{vir-change}
Let $\ls=\ls_3(d,m_1,\ldots,m_r)$ be a linear system such that $2d
\geq m_i + m_j + m_k$ for any choice of
$\{i,j,k\}\subset\{1,2,3,4\}$ then
\begin{equation}\label{vc}
v(\cre(\ls)) - v(\ls) = \sum_{t_{ij} \geq
2}\binom{1+t_{ij}}{3}-\sum_{t_{ij}\leq -2}\binom{1-t_{ij}}{3},
\end{equation}
where $t_{ij}=m_i+m_j-d$.
\end{propo}
\begin{proof}
We can reduce to the case $\ls=\ls_3(d,m_1,\ldots,m_4)$, since the
contribution of the other points is $0$ in equation~\ref{vc}. In
this case, $\dim\ls=\# T$, where $T$ is the set defined
in~\ref{T}. In order to evaluate this number, consider the sets
\begin{eqnarray*}
\Delta & := & \{(a_1,\ldots,a_4)\in\z_{\geq 0}^4\ |\ \sum_{i=1}^4
a_i=d\},
\\
{\Delta}_i & := & \{(a_1,\ldots,a_4)\in\Delta\ |\ a_i\geq d-m_i+1
\},
\end{eqnarray*}
and observe that $T = \Delta\setminus\cup_{i=1}^4\Delta_i$. By the
inclusion-exclusion principle,
\[
\#T =  \#\Delta -\sum_{i}\#\Delta_i+\sum_{i<j}\#(\Delta_i\cap
\Delta_j) - \sum_{i<j<k}\#(\Delta_i\cap \Delta_j\cap\Delta_k).
\]
The virtual dimension $v(\ls)$ is given by the first two terms on
the right side of the equation above. A point $(a_1,\ldots,a_4)$
belongs to the intersection of the first three $\Delta_i$'s if
$a_i>d-m_i$ for $i=1,2,3$; summing up these inequalities one
obtains $d-a_4>3d-m_1-m_2-m_3$ which implies that
$2d-m_1-m_2-m_3<0$, a contradiction to our assumption. The same
argument holds for any other triple of points, so this gives the
following formula
\[
\dim\ls = v(\ls) + \sum_{i<j}\#(\Delta_i\cap \Delta_j).
\]
The set $\Delta_1\cap\Delta_2$ is in one to one correspondence
with the tetrahedron
\[
\Delta_{12}:=\{(b_1,\ldots,b_4)\in\z_{\geq 0}^4\ |\ \sum_{i=1}^4
b_i=m_1+m_2-d-2\}
\]
where $b_i=a_i-d-1+m_i$ for
$i=1,2$ and $b_i=a_i$ for $i=3,4$. The same holds for the other
intersections and, since $m_i+m_j-d-2=t_{ij}-2$, we have
$\#(\Delta_i\cap\Delta_j)=\binom{t_{ij}+1}{3}$. The equality
$\dim\ls = \dim\cre(\ls)$ implies that
\[
v(\cre(\ls))-v(\ls) = \sum_{t_{ij}\geq 2}\binom{t_{ij}+1}{3} -
\sum_{t'_{ij}\geq 2}\binom{t'_{ij}+1}{3},
\]
where $t'_{ij}$ is defined in the same way of $t_{ij}$ for the
system $\cre(\ls)$. This implies that
$t'_{12}=(m_1+k)+(m_2+k)-(d+k)=t_{12}+k=-t_{34}$ and an analogous
equality holds for each $t'_{ij}$.
\end{proof}

\begin{coro}\label{decreasing-degree}
Under the same assumptions of Proposition~\ref{vir-change}, if the
degree of $\cre(\ls)$ is smaller than that of $\ls$, then
$v(\cre(\ls))\geq v(\ls)$.
\end{coro}
\begin{proof}
The difference between the degree of $\cre(\ls)$ and that of $\ls$
is equal to $k = 2d - \sum_{i=1}^4m_i$. From $2d <
\sum_{i=1}^4m_i$ we deduce that, if $t_{12}\geq 2$ then $d-m_3-m_4
< m_1+m_2-d$ which is equivalent to $-t_{34} < t_{12}$. The same
holds for each $t_{ij}$ such that $t_{ij}\geq 2$, hence the right
side of equation~\ref{vc} is non negative.
\end{proof}

\section{The isomorphisms of ${\rm A}^i\thinspace (Y)$ induced by the cubic transformation}

The resolution of the indeterminacy of map~\ref{cubic} may be
expressed by the following diagram:
\begin{equation}\label{crem}
\xymatrix{
& {}Y \ar[dl]_{p} \ar[dr]^{p'} &\\
X \ar[d]_{\pi} \ar@{-->}[rr] & {} & {} X \ar[d]^{\pi'}\\
\p^3 \ar@{-->}[rr]^{\cre} & {} &  {}\p^3}
\end{equation}
Here $\pi$ is the blow-up of $\p^3$ along the points
$p_1,\ldots,p_4$ and $p$ is the blow-up of $X$ along the strict
transforms of the lines $l_{ij}$ through $p_i$ and $p_j$. The
Picard group of $X$ is generated by $\langle
H,E_1,\ldots,E_4\rangle$ where $E_i$ is the exceptional divisor
corresponding to $p_i$. Let us denote by $F_i=p^* E_i$, this means
that $F_i$ is the blow-up of a plane through three points, and let
$F_{ij}$ be the exceptional divisor corresponding to the blow-up
of the line through $E_i$ and $E_j$. The divisor $F_{ij}$ is a
quadric and we will denote by $f_{ij}^v$ its ``vertical'' ruling,
i.e. the one given by the lines of the blow-up, and by $f_{ij}^h$
the other one (the ``horizontal''). The morphism $p'$ contracts
the quadrics $F_{ij}$ along the horizontal rulings while the
divisors contracted by $\pi'$ are the four fundamental planes. A
base for $\pic(Y)$ is given by $\langle H, F_1,\ldots,
F_4,F_{12},\ldots,F_{34}\rangle$, where, with abuse of notation,
we write $H$ instead of $p^*H$. Let $H_1 =
H-F_2-F_3-F_4-F_{23}-F_{24}-F_{34}$ be the strict transform of the
plane of $\p^3$ through the points $p_2,p_3,p_4$ and the three
lines of the triangle. In the same way we define $H_i$ for
$i=2,3,4$. The involution $\cre: \pic(Y)\rt\pic(Y)$ induced
by~\ref{cubic} is given by
\begin{equation}\label{invo1}
\cre(F_i)  =  H_i,\hspace{5mm} \cre(F_{ij})  =  F_{hk},
\end{equation}
where $\{h,k\} = \{1,2,3,4\}\setminus\{i,j\}$, which gives
\begin{equation}\label{invo2}
\cre(H)=3H-\sum_i 2F_i -\sum_{i<j} F_{ij}.
\end{equation}
\begin{propo}\label{pic}
Let $Z$ be the $0$-dimensional scheme of four points
$p_1,\ldots,p_4$ with multiplicities $m_1,\ldots,m_4$ and $W$ be
the $1$-dimensional scheme of the lines $\langle p_i,p_j\rangle$
with multiplicities $n_{ij}$. The map~\ref{cubic} transforms the
sheaf $\mathcal{F}=\oc_{\p^3}(d)\otimes\ci_Z\otimes\ci_{W}$ into
$\mathcal{F}'=\oc_{\p^3}(d+s)\otimes\ci_{Z'}\otimes\ci_{W'}$,
where $s=2d - \sum m_i$, $m'_i=m_i+s$ and
$n'_{ij}=d-m_i-m_j+n_{hk}$.
\end{propo}
\begin{proof}
Let us take the pull-back of $\mathcal{F}$ on $Y$:
\begin{eqnarray*}
p^*\pi^*\mathcal{F} & = & dH-\sum_{i} m_iF_i-\sum_{i<j}
n_{ij}F_{ij}.
\end{eqnarray*}
By~\ref{invo1} and~\ref{invo2} its image under the map $\cre$ is
\begin{eqnarray*}
\cre(p^*\pi^*\mathcal{F}) & = & (d+s)H-\sum_{i} (m_i+s)F_i -
\sum_{i<j}(d-m_i-m_j+n_{hk})F_{ij},
\end{eqnarray*}
which gives the thesis.
\end{proof}

To each linear system $\ls=\ls_3(d,m_1,\ldots,m_4)$ we associate a
$1$-cycle defined as
\[
\Gamma(\ls):=\sum_{t_{ij}\geq 1} t_{ij}l_{ij},
\]
where $t_{ij}=m_i+m_j-d$ and $l_{ij}$ is the line through $p_i$
and $p_j$. Observe that by definition
$H^0(\ls\otimes\ci_{\Gamma(\ls)}) = H^0(\ls)$, since each line
$l_{ij}\in\Gamma(\ls)$ is contained into the base locus of $\ls$
with multiplicity at least $t_{ij}$.

\begin{propo}\label{vir=}
With the preceding notations we have
\[
\chi(\ls\otimes\ci_{\Gamma(\ls)})=\chi(\ls) + \sum_{t_{ij}\geq
2}\binom{t_{ij}+1}{3}.
\]
\end{propo}
\begin{proof}
Consider the pull-back of $\ls\otimes\ci_{\Gamma(\ls)}$ on $Y$. In
what follows we will concentrate on one of the six lines and in
order to simplify the notations we will omit the indexes. Let
$l\in\ac^2(X)$ be the strict transform of the line, then
$N_{l|X}\cong\oc_{\p^1}(-1)\oplus\oc_{\p^1}(-1)$. From the
evaluation of the tautological line bundle associated to the
blow-up of $l$ (see~\cite{fu}) one obtains that ${F}_{|F} = - f^v
- f^h$. The intersection $\ls l=-t$ gives $p^*\ls_{|F} = -t f^v$.
The exact sequence
\[
\xymatrix{ 0\ar[r] & p^*\ls-(k+1)F\ar[r] & p^*\ls-kF \ar[r] &
(p^*\ls-kF)_{|F}\ar[r] & 0, }
\]
and the preceding formulas imply that
$h^0((p^*\ls-kF)_{|F})=h^0(kf^h+(k-t)f^v)=0$ if and only if $k<t$.
Therefore $tF\subseteq\bs (p^*\ls)$ and
\[
\chi(p^*\ls-tF) = \chi(p^*\ls)-\sum_{k=0}^{t}\chi(kf^h+(k-t)f^v)).
\]
By Riemann-Roch theorem on the quadric $F$ we have
$\chi(kf^h+(k-t)f^v))=(k+1)(k+1-t)$. An easy calculation shows
that the last sum of the above equation is equal to
$-\binom{t+1}{3}$, hence applying this procedure to each one of
the $F_{ij}$ one obtains:
\[
\chi(p^*\ls-\sum_{t_{ij}\geq 2}
t_{ij}F_{ij})=\chi(p^*\ls)+\sum_{t_{ij}\geq 2}
\binom{t_{ij}+1}{3}.
\]
\end{proof}

\begin{coro}\label{speciality}
Let $\ls=\ls_3(d,m_1,\ldots,m_r)$ be a linear system and
$C_1,\ldots,C_n\in\ac^2(X)$ be a set of irreducible rational
curves such that
$N_{C_i|X}\cong\oc_{\p^1}(-1)\oplus\oc_{\p^1}(-1)$ and $\ls
C_i=-t_i\leq -2$ for $i=1\ldots,n$. Then
\[
\dim\ls - v(\ls) \geq \sum_{i=1}^n\binom{t_i+1}{3} -
h^2(\ls\otimes\ci_{\Gamma}),
\]
where $\Gamma:=\sum t_iC_i$ is the $1$-dimensional scheme of the
multiple curves.
\end{coro}
\begin{proof}
Observe that in the proof of Proposition~\ref{vir=} the only
assumption needed on the curve $l$ is that its normal bundle in
$X$ is $\oc_{\p^1}(-1)\oplus\oc_{\p^1}(-1)$. Since we are assuming
the same for the $C_i$, we have that
$\chi(\ls\otimes\ci_{\Gamma})=\chi(\ls) + \sum\binom{t_i+1}{3}$.
The equalities $h^0(\ls) = h^0(\ls\otimes\ci_{\Gamma})$ and
$h^2(\ls) = 0$ give us the thesis.
\end{proof}

\begin{remark}
Let us take a system $\ls$, chose four points and consider the
associated sheaf $\ls\otimes\ci_{\Gamma(\ls)}$. Use
Proposition~\ref{pic} to transform this sheaf into
$\ls'\otimes\ci_{\Gamma(\ls')}$. The system $\ls'$ is obtained
from $\ls$ by applying Proposition~\ref{cre-surfaces}. Observe
that
$h^i(\ls\otimes\ci_{\Gamma(\ls)})=h^i(\ls'\otimes\ci_{\Gamma(\ls')})$
since the pull-back of the first system on $Y$ is just the
pull-back of the second one obtained by a base change on
$\pic(Y)$. In particular this implies that the virtual dimensions
of the two systems are the same. This, together with
Proposition~\ref{vir=} gives another proof of
Proposition~\ref{vir-change}.
\end{remark}

Let us consider now the action of transformation~\ref{cubic} on
the curves of $Y$. A basis for $\ac^2(Y)$ may be given by:
$\langle h,f_1,\ldots,f_4,f_{12},\ldots,f_{34}\rangle$, where $h$
is the pull-back of the class of a line in $\p^3$, $f_i$ is the
pull-back of a line of $E_i$ and $f_{ij}=f_{ij}^v$ is the vertical
ruling of $F_{ij}$.

\begin{propo}\label{a1a2}
The intersection matrix $M:\ac^1(Y)\times\ac^2(Y)\rt\z$ with
respect to the chosen bases of $\ac^i(Y)$ is given by
\[
M = \left(
\begin{array}{c|c}
1 & 0 \\
\hline \\
0 & -I \\
\end{array}
\right)
\]
\end{propo}
\begin{proof}
By using the projection formula, since $(p\circ\pi)_*F_i = p_i$,
$F_i$ has a non zero intersection only with $f_i$ and $F_if_i =
{F_i}_{|F_i}f_i = -f_i^2 = -1$ where the second equation is due to
$F_i^2 = p^*E_i^2 = p^* N_{E_i|X} = p^*(-e_i) = -f_i$. In the same
way, since $(p\circ\pi)_*F_{ij} = l_{ij}$, $F_{ij}$ has non
vanishing intersection only with $f_{ij}$. Let $F$ and $f$ be
$F_{ij}$ and $f_{ij}$ respectively, then $F f = F_{|F} f =
(-f^v-f^h)f^v = -1$. Finally $Hh=1$ since it is the intersection
of a plane with a line.
\end{proof}

\begin{propo}
With the same notation as before, we have that $HF_{ij}=f_{ij}$
and $F_{ij}^2=-h+f_i+f_j-2f_{ij}$.
\end{propo}
\begin{proof}
As before we will use $F,f$ instead of $F_{ij},f_{ij}$. The first
equality is obvious, since $p_*H$ intersect $p_*F$ along a point
and $f$ is just the vertical fiber of the quadric over that point.
Observe that since $F^2$ has a non vanishing intersection only
with $H,F_i,F_j$ and $F$, then by Proposition~\ref{a1a2} we have
$F^2=ah+m_if_i+m_jf_j+bf$. The coefficient $a$ is equal to $F^2H =
Ff = -1$. In the same way $m_i = - F^2F_i = - Ff = 1$ and $b = -
F^3 = - (-f^v-f^h)^2 = - 2$.
\end{proof}

\begin{propo}\label{a1a1}
The action of map~\ref{cubic} on $\ac^2(Y)$ is given by
\[
\cre(h) = 3h-\sum_{i} f_i,\ \cre(f_i) = 2h - \sum_{j\neq i} f_j,\
\cre(f_{ij}) = h + f_{ks} - f_k - f_s,
\]
where $\{k,s\} =
\{1,2,3,4\}\setminus\{i,j\}$.
\end{propo}
\begin{proof}
Since $h=H^2,\ f_1=-F_1^2$ and $f_{12}=HF_{12}$, we have
\begin{eqnarray*}
\cre(h) & = & \cre(H)^2 \\
& = & (3H-2\sum_{i} F_i-\sum_{i<j} F_{ij})^2 \\
& = & 9h-4\sum_{i} f_i -\sum_{i<j} (h-f_i-f_j+2f_{ij}) -
6\sum_{i<j} f_{ij} + 8\sum_{i<j} f_{ij} \\
& = & 3h - \sum_{i} f_i \\
& & \\
\cre(f_1) & = & -\cre(F_1)^2 \\
& = & -(H-F_2-F_3-F_4 -F_{23}-F_{24}-F_{34})^2 \\
& = & -(h-f_2-f_3-f_4 -3h +2(f_2+f_3+f_4) -
2(f_{23}+f_{24}+f_{34}) \\
& &  - 2(f_{23}+f_{24}+f_{34}) +
4(f_{23}+f_{24}+f_{34})) \\
& = & 2h-f_2-f_3-f_4. \\
& & \\
\cre(f_{12}) & = & \cre(H) \cre(F_{12}) \\
& = & (3H-2\sum_{i} F_i-\sum_{i<j} F_{ij})F_{34} \\
& = & 3f_{34}-4f_{34}+h-f_3-f_4+2f_{34} \\
& = & h+f_{34}-f_3-f_4.
\end{eqnarray*}
\end{proof}

Now, with $\ell_3(\delta,\allowbreak \mu_1,\allowbreak
\ldots,\mu_4,\beta_{12},\ldots,\beta_{34})$ we will mean the
system of curves of $\p^3$ of degree $\delta$ with multiplicities
$\mu_i$ in the four points and intersecting the line $l_{ij}$
along $\beta_{ij}$ points.

\begin{propo}\label{cre-curves-Y}
The image of the curve $\ell=\allowbreak\ell_3(\delta,\allowbreak
\mu_1,\allowbreak \ldots,\mu_4,\beta_{12},\ldots,\beta_{34})$ by
transformation~\ref{cubic} is
$\cre(\ell):=\ell_3(\delta',\mu_1',\ldots,\mu_4',\beta_{12}',\ldots,\beta_{34}')$,
where
\begin{equation}\label{crea2}
\delta' = 3\delta-\sum_{i}2\mu_i-\sum_{i<j}\beta_{ij},\hspace{5mm}
\mu_r' = \delta-\sum_{j\neq r}\mu_j-\sum_{\{i,j\}\not\ni
r}\beta_{ij},\hspace{5mm} \beta_{ij}' = \beta_{ks}.
\end{equation}

\end{propo}
\begin{proof}

The strict transform of a general element of $\ell$ may be written
on $Y$ as $\delta h - \sum\mu_i f_i - \sum\beta_{ij} f_{ij}$. The
image of this element by the map $\cre$ is
\[
\delta (3h - \sum_{i} f_i) - \sum_{i}\mu_i (2h-\sum_{j\neq i}f_j)
- \sum_{i<j}\beta_{ij} (h+f_{ks}-f_k-f_s).
\]
This gives immediately the value of $\delta'$ and $\beta_{ij}'$.
For $\mu_r'$, observe that the $f_r$ appear in the first sum with
coefficient $-\delta$, in the second with coefficient $\sum_{j\neq
r}\mu_j$ and in the third sum it appears each time that the index
$r\in\{k,s\}$ or, which is the same, each time that
$r\notin\{i,j\}$.
\end{proof}

If the system $\ell$ has no intersection with any one of the six
lines, then we have the following:

\begin{coro}\label{cre-curves}
The action of transformation~\ref{cubic} on
$\ell=\allowbreak\ell_3(\delta,\allowbreak \mu_1,\allowbreak
\ldots,\mu_r)$ is given by:
\begin{eqnarray}\label{cre-a2}
\cre(\ell) & := & \ell(\delta+2h,\mu_1+h,\ldots,
\mu_4+h,\mu_5,\ldots,\mu_r),
\end{eqnarray}
where $h=\delta-\sum_{i=1}^4 \mu_i$.
\end{coro}

The action of $\cre$ on the space of curves may be extended to a
linear action on the space of polynomials
$\co[\delta,\mu_1,\ldots,\mu_r]$. An easy calculation gives

\begin{coro}\label{invariants-2}
The following polynomials are invariant with respect to
transformation~\ref{cre-a2}:
\[
2\delta - \sum_{i=1}^s\mu_i,\hspace{5mm} \delta^2 -
\sum_{i=1}^s2\mu_i^2 + 3.
\]
\end{coro}

\section{Conjecture}
In this section we state a conjecture which allows us to give a
procedure for computing the dimension of a linear system
$\ls=\ls_3(d,m_1,\ldots,m_r)$. First of all, since birational
transformations do not change the effective dimension of $\ls$, we
can perform Cremona transformations until we get a system in
standard form. If at some step we get a system $\ls^k$ with a
negative multiplicity $-\alpha_i$ at $p_i$, then the system
contains $\alpha_iE_i$ in its base locus. We remark that if this
is the case, then there exists also a fixed component contained
$\alpha_i$ times in the base locus of the starting system $\ls$.
This is the image of $E_i$ by the sequence of Cremona
transformations sending $\ls^k$ back to $\ls$. This component can
be removed without changing the dimension of $\ls$. In particular
We can remove $\alpha_iE_i$ from $\ls^k$ and keep performing
Cremona if possible. We can then reduce to the study of the
dimension of linear systems in standard form. In this direction
let us state the following:

\begin{con}\label{conj}
A linear system $\ls=\ls_3(d,m_1,\ldots,m_r)$ in standard form is
special if and only if one of the following holds:
\begin{itemize}
\item[(i)] there exists a quadric $Q=\ls_3(2,1^9)$ such that $Q(\ls-Q)(\ls-K)<0$;
\item[(ii)] there exists a line $\ell=\ell_3(1,1^2)$ such that $\ls\ell \leq -2$.
\end{itemize}
\end{con}
We remark that if condition (i) holds, then by equation~\ref{add}
$v(\ls)<v(\ls-Q)$, while $\dim\ls\geq\dim(\ls-Q)$, which means
that $\ls$ is special. In order to simplify the procedure we are
going to prove that, under an extra assumption, condition (i)
implies that the quadric $Q$ is contained in the base locus of
$\ls$.

\begin{lemma}\label{HHQ}
Let us suppose that Harbourne-Hirschowiz Conjecture holds for
linear systems on $\p^2$ with $10$ fixed points. If
$\ls=\ls_3(d,m_1,\ldots,m_r)$ is a system in standard form and
such that $Q(\ls-Q)(\ls-K) < 0$, then $Q\subset\bs(\ls)$.
\end{lemma}
\begin{proof}
From the exact sequence
\[
\xymatrix{ 0\ar[r] & \ls-Q \ar[r] & \ls \ar[r] & \ls_{|Q}\ar[r] &
0, }
\]
we get $v(\ls)=v(\ls-Q)+v(\ls_{|Q})+1$. If we compare with
equation~\ref{add} we obtain that
\[
v(\ls_{|Q}) = \frac{Q(\ls-Q)(\ls-K)}{2} - 1 < 0.
\]
The system $\ls_{|Q}=\ls_Q((d,d),m_1,\ldots,m_9)$ is equivalent to
the planar system $\ls_2=\ls_2(2d-m_1,(d-m_1)^2,m_2,\ldots,m_9)$
(see~\cite{au2}). Therefore in order to prove the thesis it is
enough to prove that the system $\ls_2$ is non-special or, since
we assume that Harbourne-Hirschowiz conjecture holds for $10$
multiple points, to prove that $\ls_2$ is not $(-1)$-special.
Let us compare the multiplicity $d-m_1$ with respect to the $m_i$'s.\\
If $d-m_1\geq m_3$, the system $\ls_2$ is Cremona-stable (since
$2d-m_1\geq 2(d-m_1)+m_2$)
and hence it is not special.\\
If $d-m_1 < m_4$ and $\ls_2$ is not Cremona-stable, then $2d-m_1 <
m_2+m_3+m_4$, which
is not possible since we are assuming that $\ls$ is in standard form.\\
Finally, if $m_4\leq d-m_1 < m_3$ and $\ls_2$ is not
Cremona-stable, then $d\leq m_2+m_3-1$. If we write $d=m_2+m_3-t$,
with $t\geq 1$, then $\ls_2=\ls_2(2m_2+2m_3-m_1-2t,m_2,m_3,
(m_2+m_3-m_1-t)^2,m_4,\ldots,m_9)$. Performing a Cremona
transformation with the first three points we obtain a stable
system (the degree equals the sum of the highest multiplicities).
\end{proof}

Now we wonder: which system of quadrics may be contained into the
base locus of a given linear system? The answer is given by the
following:
\begin{lemma}\label{baseq}
If $|\sum_{i=1}^nr_iQ_i|$ is contained in the base locus of a
linear system, then the $Q_i$ must share $8$ points, i.e.
$Q_i\in\ls_3(2,1^8)$ for each $i$.
\end{lemma}
\begin{proof}
Consider any two of the given quadrics, say $Q_1,Q_2$, which share
$s\leq 8$ points. Since $|Q_1+Q_2|=\ls_3(4s,2^s,1^{18-2s})$ has
virtual dimension $16-2s$, the system moves unless $s=8$. Now
consider any three quadrics, say $Q_1,Q_2,Q_3$. If they have less
than $8$ common points, then $|Q_1+Q_2+Q_3|=\ls_3(6,3^7,2^3)$
since, by the preceding discussion, each pair must share $8$
points. But this system has virtual dimension $1$, which is not
possible. If four or more quadrics share less than $8$ points,
then there exist three of them which do not share $8$ points, and
this is not possible as proved before. \\
On the other hand, let us prove that if $s=8$, then the system
$\ls=|\sum_{i=1}^nr_iQ_i|$ cannot move. We can write
$\ls=\ls_3(2r,r^8,r_1,\ldots,r_n)$, where $r=\sum_{i=1}^nr_i$.
From the exact sequence of $Q_1$ we get
\[
\xymatrix{ 0\ar[r] & \ls-Q_1 \ar[r] & \ls \ar[r] &
\ls_{|Q_1}\ar[r] & 0. }
\]
Let us consider the restricted system $\ls_{|Q_1}=\ls_Q((2r,2r),
r^8,r_1)$. By~\cite{au2} this is equivalent to the planar system
$\ls_2(3r,r^9,r_1)$, which is empty since $\ls_2(3r,r^9)$ is the
fixed cubic $\ls_2(3,1^9)$ counted $r$ times. Therefore also
$\ls_{|Q_1}=\emptyset$ which implies that $Q_1$ is contained in
the base locus of $\ls$. This also implies that $r_1Q_1$ is
contained in the base locus of $\ls$ and the same holds for
$r_iQ_i$. Therefore $\ls = \sum_{i=1}^nr_iQ_i$.
\end{proof}
In order to conclude the part of the procedure concerning quadrics
we need the following:
\begin{lemma}
Let us suppose that Harbourne-Hirschowiz conjecture holds for
linear systems on $\p^2$ with $10$ points. If
$\ls=\ls_3(d,m_1,\ldots,m_r)$ is in standard form and
$Q(\ls-Q)(\ls-K) < 0$, where $Q$ is the quadric through the first
$9$ points, then $\ls-Q$ is still in standard form.
\end{lemma}
\begin{proof}
By hypothesis $2d\geq\sum_{i=1}^4 m_i$; moreover the degree of
$\ls-Q$ is $d-2$ while the first $9$ multiplicities are $m_i-1$.
Therefore $\ls-Q$ is in standard form unless $m_{10} > m_4-1$. But
this can happen only if $m_i=m$ for $i=4,\ldots 10$. In this case
$\ls$ would contain all the quadrics through $p_1,p_2,p_3$ and $6$
of the points $p_4,\ldots,p_{10}$, but by Lemma~\ref{baseq} this
is not possible.
\end{proof}
Therefore, assuming that Harbourne-Hirschowiz conjecture holds up
to $10$ points, we can proceed computing $Q(\ls-Q)(\ls-K)$. If it
is negative we can remove the quadric $Q$ and consider the system
$\ls-Q$, which is still in standard form. We reorder the
multiplicities if necessary and we keep removing the quadric as
far as $Q(\ls-Q)(\ls-K) < 0$.

Let us consider part (ii) of Conjecture~\ref{conj}. We are going
to prove that, in fact, such a system $\ls$ is special.

\begin{propo}\label{h2=0}
Let $\ls=\ls_3(d,m_1,\ldots,m_r)$ be a non-empty linear system and
let $l$ be the line through $p_1,p_2$. If $\ls l=-t\leq -1$ then
$h^2(\ls\otimes\I_{tl})=0$.
\end{propo}
\begin{proof}
Let $Z=Z'+Z''$, where $Z'=m_1p_1+m_2p_2$ and
$Z''=m_3p_3+\ldots+m_rp_r$ and let $\ls'=
\oc_{\p^3}(d)\otimes\I_{Z'}$. The tensor product of the defining
sequence of $Z''$ with $\ls'\otimes\I_{tl}$ gives
\[
\xymatrix{ 0\ar[r] & \ls\otimes\I_{tl} \ar[r] & \ls'\otimes\I_{tl}
\ar[r] & \oc_{Z''}\ar[r] & 0. }
\]
Since $h^i(\oc_{Z''})=0$ for $i\geq 1$, it is enough to prove that
$h^2(\ls'\otimes\I_{tl})=0$. \\
Let us take a plane $V\in\ls(1,1^2)$ and denote by $W$ the
corresponding element of $\ls(1,1^2)\otimes\ci_l$. From the
defining sequence of $W$ we obtain
\[
\xymatrix{ 0\ar[r] & (\ls'-V)\otimes\I_{(t-1)l} \ar[r] &
\ls'\otimes\I_{tl} \ar[r] & (\ls'\otimes\I_{tl})_{|W} \ar[r] & 0.
}
\]
Observe that
$h^2((\ls'\otimes\I_{tl})_{|W})=h^2(\oc_{\p^2}(d-t)\otimes\ci_{Z'})=0$,
since $t=m_1+m_2-d\leq d$ (otherwise $m_1 > d$ and the system
would be empty). This means that
$$h^2((\ls'-V)\otimes\I_{(t-1)l})\geq h^2(\ls'\otimes\I_{tl}),$$
so we can proceed by induction on $t$ until we obtain the system
$\ls_3(d-t,m_1-t,m_2-t)$ whose $h^2$ vanishes.
\end{proof}
An easy consequence of Propositions~\ref{speciality}
and~\ref{h2=0} is the following:
\begin{coro}\label{line}
Let $\ls$ and $l$ be as before. If $\ls l=-t\leq -2$, then $\ls$
is special and $\dim\ls - v(\ls)\geq\binom{t+1}{3}$.
\end{coro}
The preceding discussion allows us to give a procedure to
calculate $\dim\ls$.
\begin{center}
\begin{table}[h]
\caption{\ }\label{pro} \fbox{
\begin{minipage}{10cm}
\vspace{2mm} \texttt{\noindent While $\ls$ is not in standard form
\begin{itemize}
\item[] put $\ls:=\cre(\ls)$.
\item[] If $m_i=-\alpha_i<0$ then put $\ls:=\ls-\alpha_iE_i$.
\end{itemize}
} \texttt{\noindent While $Q(\ls-Q)(\ls-K)< 0$
\begin{itemize}
\item[] put $\ls:=\ls-Q$.
\end{itemize}
} \texttt{\noindent Put $t_{ij}:=m_i+m_j-d$. Return $
\displaystyle\dim\ls=v(\ls)+\sum_{t_{ij}\geq
2}\binom{t_{ij}+1}{3}$.} \\
\end{minipage}
}
\end{table}
\end{center}
Observe that in the last step of the procedure we assume that
$h^2(\ls\otimes\ci_{\Gamma})=0$, where $\Gamma=\sum t_{ij}l_{ij}$
is the $1$-dimensional scheme of the lines $l_{ij}$, such that
$\ls l_{ij}=-t_{ij}<0$.

\section{Motivations}

In~\cite{dl} Conjecture~\ref{conj} has been proved for systems
$\ls_3(d,m_1,\ldots,m_r)$ with $r\leq 8$. With the help of a
procedure written with Singular~\cite{sing}, we verified the
conjecture for systems $\ls_3(d,m^r)$ with $m\leq 7$
and $r\leq 20$. \\
A consequence of the conjecture is that the only curves $\ell$
which can give speciality are the images of the line
$\ell_3(1,1^2)$ by a finite set of Cremona transformations.  In
one direction, observe that if there exists a curve $\ell$ which
is the image of $\ell_3(1,1^2)$ by a finite set of Cremona
transformations and such that $\ls\ell=-t\leq -2$, then $\ls$ is
special. In fact, performing back Cremona, we get a system $\ls'$
such that $\ls'\ell_3(1,1^2)=-t$, since the intersection product
is invariant under Cremona. Therefore, by Proposition~\ref{h2=0},
since cohomology groups are invariant under birational
transformations, $h^2(\ls\otimes\I_{t\ell})=0$, and $\ls$ is
special. On the other hand, for other classes of curves one can
have the following problem:
\begin{example}
Let $\ell=\ell_3(4,1^8)$ be the system of quartics through $8$
fixed points and suppose that $\ls\ell=-t\leq -2$. The
intersection product is $4d-\sum_{i=1}^8m_i\leq -2$, and hence it
must be $2d-\sum_{i=1}^4 m_i\leq -1$, which implies that $\ls$ is
not in standard form. Let us perform a Cremona transformation
involving the first $4$ points. Since the intersection product is
invariant under this transformation, $\ls'\ell'=\ls\ell=-t$, where
$\ls'=\cre(\ls)$ and $\ell'=\cre(\ell)$. Moreover the degree $d'$
of $\ls'$ is strictly smaller than $d$ and, by
Proposition~\ref{cre-curves} $\ell'=\ell$. Arguing as before we
can say that $\ls'$ is not in standard form and we can perform
another transformation. Iterating this procedure we can transform
$\ls$ until we get a system having negative degree, and hence we
conclude that it must be $\ls=\emptyset$. Therefore the quartic
can not give speciality. The fact that $\ell_3(4,1^8)$ can not be
obtained from $\ell_3(1,1^2)$ by a finite set of Cremona
transformations is an easy consequence of
Corollary~\ref{invariants-2}, since the invariant $\delta^2-\sum
2\mu_i^2 +3$ is $3$ for the quartic and $0$ for the line.
\end{example}
Another consequence of Conjecture~\ref{conj} is that if $\ls$ is
in standard form then we can not find any curve
$\ell\neq\ell_3(1,1^2)$ which is the image of a line and such that
$\ls\ell\leq -2$. In order to give a motivation for this we prove
the following:
\begin{lemma}\label{inter}
Let us suppose that $\ls=\ls_3(d,m_1,\ldots,m_r)$ and $\ell$ is
obtained from $\ell_3(1,1^2)$ by a finite set of Cremona
transformations such that at each step the degree of the curve
increases. Then we can write
\begin{equation}\label{int}
\ls\ell=\beta_1(2d-\sum_{i=1}^4m_i^{(1)})+\cdots+
\beta_a(2d-\sum_{i=1}^4m_i^{(a)})+(d-m_h-m_k),
\end{equation}
where $\beta_j\geq 1$, $m_i^{(j)}$ are chosen between $m_1,\ldots,
m_r$ and $h,k\geq 5$.
\end{lemma}
\begin{proof}
We argue by induction on the number of Cremona transformations
necessary to obtain $\ell$ from the line $\ell_3(1,1^2)$. First of
all, after one transformation the image of the line is the
rational normal cubic $\ell_3(3,1^6)$, having intersection product
with $\ls$ equal to $3d-\sum_{i=1}^6 m_i=(2d-\sum_{i=1}^4
m_i)+(d-m_5-m_6)$. Now we assume that the formula is true for
$\ell=\ell_3(\delta,\mu_1,\ldots,\mu_s)$ and we prove it for the
curve $\ell'=\cre(\ell)=\ell_3(\delta',\mu_1',\ldots,\mu_s')$
obtained from $\ell$ performing one more Cremona transformation
increasing the degree. We can suppose that the transformation is
based on the first $4$ points. By formula~\eqref{cre-a2},
$\delta'=\delta + 2h$ and $\mu_i'=\mu_i + h$ for $i=1,\ldots,4$,
where $h=\delta-\sum_{i=1}^4 \mu_i>0$, and $\mu_i'=\mu_i$ for
$i\geq 5$. Therefore $\ls\ell'-\ls\ell=h(2d-\sum_{i=1}^4 m_i)$,
which gives the thesis.
\end{proof}
\begin{coro}
If $\ls$ is not empty and in standard form and $\ell$ is obtained
from $\ell_3(1,1^2)$ by a finite set of Cremona transformations
such that at each step the degree of the curve increases, then
$\ls\ell\geq 0$.
\end{coro}
\begin{proof}
Since $\ls$ is in standard form, $2d\geq\sum m_i^{(j)}$. Moreover
$d-m_h-m_k\geq 0$ since otherwise $2d < m_1 + m_2 + m_h + m_k$.
Hence all the terms on the right side of equation~\ref{int} are
non-negative.
\end{proof}
Therefore, if we assume that the following conjecture holds,
we have that a system in standard form can not
have negative intersection product with the image of a line.
\begin{con}\label{curves}
Let $\ell=\ell_3(\delta,\mu_1,\ldots,\mu_r)$ be a curve that can
be obtained from the line $\ell_3(1,1^2)$ by a finite set of
Cremona transformations. Then $\ell$ can be obtained by a finite
set of Cremona transformations such that at each step the degree
of the curve increases.
\end{con}

\section{Homogeneous Linear Systems}

In this section we study special homogeneous systems
$\ls=\ls_3(d,m^r)$.
\begin{propo}
The system $\ls$ is empty for $d\leq 2m-1$ and $r\geq 8$.
\end{propo}
\begin{proof}
It is enough to show that the system
$\ls_3(2m-1,m^8)$ is empty. We first prove by induction that
performing Cremona on the first four points and the last four
alternatively, at the $i$-th step we get the system
$\ls^i=\ls_3(2m-2i^2-1,(m-i^2+i)^4,(m-i^2-i)^4)$. The basis of
induction holds since applying Cremona to the starting system we
get $\ls^1=\ls_3(2m-3,(m-2)^4,m^4)$. If we now perform Cremona to
$\ls^i$ taking the first four points we have that
$k=2(2m-2i^2-1)-4(m-i^2+i)=-4i-2$ and hence the transform system
is
$\cre(\ls^i)=\ls_3(2m-2(i+1)^2-1,(m-(i+1)^2-(i+1))^4,(m-i^2-i)^4)$.
The multiplicity $m-i^2-i$ can be written as $m-(i+1)^2+(i+1)$ and
hence, reordering the points,
we get that $\cre(\ls^i)=\ls^{i+1}$, which gives the inductive step.\\
We keep performing these transformations until $m-i^2-i\leq 0$. At
this step we have that $2m-2i^2-1 < m-i^2+i$ and hence the system
$\ls^i$ is empty since it has some multiplicity bigger than the
degree.
\end{proof}
By assuming Conjecture~\ref{conj} and Harbourne-Hirschowiz
conjecture for linear systems on $\p^2$ with $10$ points, we can
prove the following:
\begin{propo}
If $d\geq 2m$ the system $\ls=\ls_3(d,m^r)$ is special if and only
if $r=9$ and $2m\leq d<[-1+\frac{3}{2}\sqrt{2m^2+2m}]$.
\end{propo}
\begin{proof}
Our assumption on the degree implies that $\ls$ is in standard
form. According to Conjecture~\ref{conj}, the system $\ls$ is
special if and only if there exists either a quadric
$Q\in\ls_3(2,1^9)$ such that $Q(\ls-Q)(\ls-K)<0$, or a line
$\ell=\ell_3(1,1^2)$ such that $\ell\ls\leq -2$. But under our
hypothesis on the degree, a line through $2$ fixed points has
non negative intersection product with $\ls$, and hence $\ls$ is special
if and only if $Q(\ls-Q)(\ls-K) < 0$. But we can write
$Q(\ls-Q)(\ls-K)=2(d-2)(d+4)-9(m-1)(m+2)$, which
is negative if and only if $2m\leq d<[-1+\frac{3}{2}\sqrt{2m^2+2m}]$.
We end the proof showing that if $\ls$ is special then the number $r$ of fixed
points can not be bigger than $9$. If we suppose by contradiction that
$r\geq 10$, then $Q_i(\ls-Q_i)(\ls-K) < 0$, where $Q_i$ is the quadric through
all the first $10$ points but $p_i$. By Lemma~\ref{HHQ}
$\ls$ must contain the system of quadrics $\sum_{i=1}^{10}Q_i$,
and this is not possible because of Lemma~\ref{baseq}.
\end{proof}

We can conclude that if the system $\ls$ has more than $9$ fixed
points (or exactly $8$ points) then it is not special. If it has
$9$ fixed points, it is special if and only if $2m\leq d<
[-1+\frac{3}{2}\sqrt{2m^2+2m}]$. If $r\leq 7$ and $d\geq 2m$, the
system can not be special. Finally, if $r\leq 7$ and $d\leq 2m-1$,
we have to follow the procedure of table~\ref{pro}.

\section{Examples}
We end the paper by presenting many examples of special systems
and computing their dimension following the procedure of
table~\ref{pro}. We remark that for each example the dimension we
found in this way agrees with the dimension computed with the help
of a Singular program.
\begin{notation}
We will use the symbols $\rri$, $\rrii$ and $\rriii$ to signify that we
are applying the first, the second and the third step of the procedure
(in (i) we are applying a Cremona transformation, in (ii)
we are removing a plane and in (iii) we are removing a quadric).
\end{notation}

\begin{example}\label{1i}
Consider the system $\ls:=\ls_3(7,4^6)$ with $v(\ls)=-1$.
Following the procedure of table~\ref{pro} we can apply four
Cremona transformations:
\[
\ls_3(7,4^6)\rri\ls_3(5,4^2,2^4)\rri\ls_3(3,2^4)\rri\ls_3(1),
\]
so that $\dim\ls = 3$. Observe that since
$\ls_3(7,4^6)\ell_3(3,1^6)=-3$ the starting system contains the
rational normal curve. After the first transformation
$\ell_3(3,1^6)\rr\ell_3(1,1^2)$ and, in fact,
$\ls_3(5,4^2,2^4)\ell_3(1,1^2)=-3$.
\end{example}

\begin{example}\label{1}
Consider the system $\ls:=\ls_3(12,7^6)$ with $v(\ls)=-50$. As
before we can apply the following Cremona transformations:
\[
\ls_3(12,7^6)\rri\ls_3(8,7^2,3^4)\rri\ls_3(4,3^4,-1^2)\rrii\ls_3(4,3^4)\rri\ls_3(0,-1^4)\rrii
\ls_3(0),
\]
so that $\dim\ls = 0$. From the procedure we may deduce that $\ls$
is given by the union of six surfaces of type $\ls_3(2,2,1^5)$ (each surface
can be obtained from a plane $\ls_3(1,1^3)$ applying one Cremona transformation).
\end{example}
In these two first examples, performing some Cremona transformations and removing
fixed planes we get a system in standard form which is not special. Let us
give two examples where the system in standard form is still special.
\begin{example}\label{13}
Consider the system $\ls:=\ls_3(10,6^5)$ with $v(\ls) = 5$. Apply
the Cremona transformation:
\[
\ls_3(10,6^5)\rri\ls_3(6,6,2^4).
\]
The last system is in standard form and has virtual dimension $11$
while its dimension is $15$, as expected by our conjecture. In fact each
of the four lines $l_{1j}$, $j=2,\ldots,5$ has intersection product $-2$
with the system.
\end{example}

\begin{example}\label{2}
Consider the system in standard form $\ls:=\ls_3(16,11,7^8)$ with
$v(\ls) = 10$. Since $Q(\ls-Q)(\ls-K) = -2$ following the
procedure, we may remove the quadric from the base locus of $\ls$:
\[
\ls_3(16,11,7^8)\rriii\ls_3(14,10,6^8).
\]
The virtual dimension of the last system is $11$ and its dimension
is $19$ (each of the $8$ lines $l_{1j}$, $j=2,\ldots,9$ has intersection
product $-2$ with the system).
\end{example}
We remark that in the two examples above, when the linear system is
in standard and it does not contain any quadric, then its speciality
is given exactly by the sum of the binomials $\binom{t_{ij}+1}{3}$,
where $-t_{ij}$ is the intersection product of the line $l_{ij}$ with the
system. In the following example we show that if a linear system $\ls$ is
not in standard form, then its speciality can be strictly smaller than
the sum of these binomials, or equivalently
$h^2(\ls\otimes\I_{\Gamma})\neq 0$ (where $\Gamma$ is the $1$-dimensional
scheme of the multiple lines).
\begin{example}\label{h2non0}
The system $\ls=\ls_3(3,3^3)$ with $v(\ls)=-11$ is not in standard
form, hence
\[
\ls_3(3,3^3)\rri\ls_3(0,-3)\rrii\ls_3(0).
\]
This implies that $\dim\ls=0$, since it is three times the plane
$\ls_3(1,1^3)$. For each line $l_{ij}$ through two of the three
points we have $\ls l_{ij}=-3$. Let
$\Gamma=3(l_{12}+l_{13}+l_{23})$. By Proposition~\ref{speciality}
$\dim\ls - v(\ls)\geq 12-h^2(\ls\otimes\I_{\Gamma})$ which
implies that $h^2(\ls\otimes\I_{\Gamma})\geq 1$. \\
\end{example}
In the next example we give a way to construct a class of
systems whose speciality is due to multiple quadrics in the fixed
locus.
\begin{example}\label{quadrics}
Let $r_1,\ldots,r_n$ be positive integers and let $r=\sum r_i$. As
we proved in Lemma~\ref{baseq}, the system
$\ls=\ls_3(2r,r^8,r_1,\ldots, r_n)$ has dimension $0$, while its
virtual dimension is $\sum( r_i - \binom{r_i+2}{3})\leq 0$.
Moreover this sum is $0$ if and only if each $r_i=1$, otherwise
the system is special.
\end{example}


\begin{thebibliography}{1}

\bibitem{ci}
Ciro Ciliberto, in {\it European Congress of Mathematics, Vol. I
(Barcelona, 2000)}, 289--316, Progr. Math., 201, Birkh\"auser,
Basel, 2001.

\bibitem{cm}
Ciro Ciliberto and Rick Miranda, \emph{Degenerations of planar
linear systems},
  J. Reine Angew. Math. \textbf{501} (1998), 191--220.

\bibitem{cm2}
\bysame, \emph{{Linear systems of plane curves with base points of
equal
  multiplicity.}}, Trans. Am. Math. Soc. \textbf{352} (2000), no.~9, 4037--4050
  (English).

\bibitem{dl}
Cindy De Volder and Antonio Laface, \emph{On linear systems of
$\p^3$ through multiple points}, preprint.

\bibitem{fu}
William Fulton, \emph{Intersection theory}, Second edition,
Springer, Berlin, 1998.

\bibitem{sing}
G.-M. Greuel, G.~Pfister, and H.~Sch\"onemann, \emph{{\sc
Singular} 2.0}, {A
  Computer Algebra System for Polynomial Computations}, Centre for Computer
  Algebra, University of Kaiserslautern, 2001, {\tt
  http://www.singular.uni-kl.de}.

\bibitem{ha}
Robin Hartshorne, \emph{Algebraic geometry}, Springer, New York,
1977; MR {\bf 57}.


\bibitem{hi}
Andr{\'e} Hirschowitz, \emph{Une conjecture pour la cohomologie
des diviseurs
  sur les surfaces rationnelles g\'en\'eriques}, J. Reine Angew. Math.
  \textbf{397} (1989), 208--213.

\bibitem{au2}
Antonio Laface and Luca Ugaglia, \emph{A counterexample to a conjecture
on linear systems on $\p^3$}, to appear on Advances in Geometry.

\bibitem{ev}
Evain Laurent, \emph{La fonction de {H}ilbert de la r\'eunion de
$4\sp h$ gros
  points g\'en\'eriques de ${\bf {p}}\sp 2$ de m\^eme multiplicit\'e}, J.
  Algebraic Geom. \textbf{8} (1999), no.~4, 787--796.

\bibitem{mi}
Thierry Mignon, \emph{Syst\`emes de courbes planes \`a
singularit\'es
  impos\'ees: le cas des multiplicit\'es inf\'erieures ou \'egales \`a quatre},
  J. Pure Appl. Algebra \textbf{151} (2000), no.~2, 173--195.



\end{thebibliography}
\end{document}